\documentclass[reqno,a4paper,12pt]{amsart} 

\usepackage{amsmath,amscd,amsfonts,amssymb}
\usepackage{mathrsfs,dsfont}

\numberwithin{equation}{section}
\numberwithin{figure}{section}

\newcommand\R{\mathbb{R}}

\newcommand\Z{\mathbb{Z}}

\newcommand\T{\mathbb{T}}

\newcommand\del{\delta}

\newcommand\lam{\lambda}
\newcommand\Lam{\Lambda}

\newcommand\sig{\sigma}

\newcommand\eps{\varepsilon}

\renewcommand\le{\leqslant}
\renewcommand\ge{\geqslant}
\renewcommand\leq{\leqslant}
\renewcommand\geq{\geqslant}
\newcommand\sbt{\subset}

\newcommand{\dotprod}[2]{\langle #1 , #2 \rangle}

\newcommand{\sign}{\operatorname{sign}}

\theoremstyle{plain}
\newtheorem{thm}{Theorem}[section]

\newtheorem{lemma}[thm]{Lemma}

\newtheorem{prop}[thm]{Proposition}

\newtheorem*{claim*}{Claim}

\newcommand{\thmref}[1]{Theorem~\ref{#1}}
\newcommand{\secref}[1]{Section~\ref{#1}}

\newcommand{\lemref}[1]{Lemma~\ref{#1}}

\theoremstyle{definition}

\newtheorem*{definition*}{Definition}
\newtheorem*{remarks*}{Remarks}
\newtheorem*{remark*}{Remark}

\begin{document}

\title[No unconditional Schauder frames of translates]
{There are no unconditional Schauder frames of translates in $L^p(\mathbb{R})$, $1\le p\le 2$}

\author{Nir Lev}
\address{Department of Mathematics, Bar-Ilan University, Ramat-Gan 5290002, Israel}
\email{levnir@math.biu.ac.il}

\author{Anton Tselishchev}
\address{Department of Mathematics, Bar-Ilan University, Ramat-Gan 5290002, Israel}
\address{St. Petersburg Department, Steklov Math. Institute, Fontanka 27, St. Petersburg 191023 Russia}
\email{celis-anton@yandex.ru}

\date{November 14, 2024}
\thanks{Research supported by ISF Grant No.\ 1044/21.}

\begin{abstract}
It is known that a system formed by translates of a single function cannot be an unconditional Schauder basis in the space $L^p(\mathbb{R})$ for any $1\le p < \infty$. To the contrary, there do exist unconditional Schauder frames of translates in $L^p(\mathbb{R})$ for every $p>2	$. The existence of such a system for $1 < p \leq 2$, however, has remained an open problem. In this paper the problem is solved in the negative: we prove that none of the spaces $L^p(\mathbb{R})$, $1\le p \le 2$, admits an unconditional Schauder frame of translates. 
\end{abstract}

\maketitle


\section{Introduction}
\label{secA1}

\subsection{}
A system of vectors $\{x_n\}_{n=1}^{\infty}$ in a Banach space $X$
is called a \emph{Schauder basis} if every  $x \in X$  admits a unique series expansion
$x = \sum_{n=1}^{\infty} c_n x_n$ where $\{c_n\}$ are scalars.
In this case there exist 
 biorthogonal continuous linear functionals $\{x_n^*\}$
such that  the coefficients of the series expansion are given by
$c_n = x_n^*(x)$   (see e.g.\ \cite[Chapter II.B]{Wojt}).
 If the series converges unconditionally
 (i.e.\ if it converges for any rearrangement of its terms)
for every  $x \in X$, then $\{x_n\}$ 
is said to be an \emph{unconditional} Schauder basis.

Given a function $g \in L^p(\R)$, we denote its translates by
\begin{equation}
\label{eqTL}
	(T_\lambda g)(x) = g(x-\lambda), \quad \lam \in \R.
\end{equation}
A long-standing problem which is still open, asks
whether the space $L^p(\R)$, $1 < p < \infty$, admits
a Schauder basis consisting of translates of a single function.
This problem was first posed for $L^2(\R)$ in \cite{OlsZal},
and later for $L^p(\R)$ in \cite[Problem 4.4]{OSSZ}.
We note that in $L^1(\R)$  there are no Schauder bases of translates,
see e.g.\ \cite[Theorem 1.7, Proposition 1.8]{OSSZ}.

It  is known that a system formed by translates of a single function 
cannot form an \emph{unconditional} Schauder basis in 
any of the spaces $L^p(\mathbb{R})$, $1< p < \infty$. 
This was proved in \cite{OlsZal} for $p=2$, in \cite{OSSZ} for $1 < p \leq 4$,
and in \cite{FOSZ} for $p>4$.

\subsection{}
If $X$ is a Banach space with dual space $X^*$, then a system
$\{(x_n, x_n^*)\}_{n=1}^\infty$ in $X \times X^*$ is called a 
\emph{Schauder frame} (or a quasi-basis) if every $x\in X$ has a series expansion
	\begin{equation}
		\label{def_Sch_frame}
	x=\sum_{n=1}^\infty x_n^* (x) x_n.
	\end{equation}
If the series \eqref{def_Sch_frame}
converges unconditionally for every  $x \in X$, then 
$\{(x_n, x_n^*)\}$ is called an \emph{unconditional} Schauder frame.
We note that if $\{x_n\}$ is a Schauder basis
with  biorthogonal coefficient functionals $\{x_n^*\}$,
then  $\{(x_n, x_n^*)\}$ is a Schauder frame. 
To the contrary, for a Schauder frame $\{(x_n, x_n^*)\}$  the series expansion 
\eqref{def_Sch_frame} need not be unique  and the coefficient functionals $\{x_n^*\}$ 
 need not be biorthogonal to $\{x_n\}$.	Hence 
	Schauder frames form a wider class of representation systems than Schauder bases.

It was proved in \cite{FOSZ} that in the space $L^p(\R)$, $p>2$,
there   exist  unconditional Schauder frames consisting 
of translates of a single function, i.e.\ of the form 
$\{(T_{\lam_n} g, g_n^*)\}$ where $g \in L^p(\R)$,
$\{\lam_n\}$ are real numbers, and
$\{g_n^*\}$ are continuous linear functionals on $L^p(\R)$.
Moreover, $\{\lam_n\}$ may be chosen to be 
\emph{an arbitrary unbounded sequence},
 and in particular, it may increase arbitrarily fast.
 
On the other hand, it has been known that in $L^1(\R)$ there are no
unconditional Schauder frames of translates, see \cite[Section 4.3]{BerCar}.

The existence of an unconditional Schauder frame of translates in $L^p(\R)$,
 $1 < p \leq 2$, has remained an open problem. 
The main goal of the present paper is to solve this problem in the negative.
We will prove the following result:

	\begin{thm}
	\label{thmA1}
		For   $1 \le p\le 2$ there does not exist 
in the space $L^p(\R)$ any unconditional Schauder frame
consisting of translates, that is,  there is no  unconditional Schauder frame of the form
$\{(T_{\lam_n} g, g_n^*)\}$  where $g \in L^p(\R)$,
$\{\lam_n\}$ are real numbers, and
$\{g_n^*\}$ are continuous linear functionals on $L^p(\R)$.
	\end{thm}

A special case of this result, where
 the coefficient functionals $\{g_n^*\}$ are assumed to be
\emph{seminormalized}, was proved in \cite{BerCar}.
 (Recall that a system of vectors $\{x_n\}$
in a Banach space $X$  is said to be  seminormalized
if there exist positive constants $A,B$ such that $A \le \|x_n\|\le B$
for every $n$.)

We note that contrary to \thmref{thmA1}, there do exist (not unconditional) Schauder frames
of translates in $L^p(\R)$ for any $1\le p < \infty$, see \cite[Section 4]{FPT}.

The rest of the paper is organized as follows. In \secref{secP1}
we give some preliminary background. In \secref{secP2} we prove
\thmref{thmA1}. Finally, \secref{secR1} contains some additional remarks.


	\section{Preliminaries}
	\label{secP1}
	
	In this section we state some known facts that will be used in the proof of \thmref{thmA1}.
	 
	 \subsection{}
	 We start with a basic property of unconditional Schauder frames, that can be proved
	 using the uniform boundedness principle in the same way 
	 as a similar statement for unconditional Schauder bases (see e.g.\ \cite[Proposition 3.1.3]{AlKal}).

	\begin{lemma}
		 	\label{unc}
			Let $\{(x_j, x_j^*)\}_{j=1}^\infty$ be an unconditional Schauder frame in 
		 a Banach space $X$.  Then there exists a constant $K$ such that for any 
		 $x\in X$ and for any
		 choice of scalars $\theta_j$ satisfying  $|\theta_j|\le 1$, we have
		 	\begin{equation}
		 	\label{eqtheta}
		\Big\| \sum_{j=1}^\infty \theta_j x_j^*(x) x_j \Big\| \leq K \|x\|
	\end{equation}
	and the series in \eqref{eqtheta} converges unconditionally.
	\end{lemma}

	We remark that the condition \eqref{eqtheta} is often 
	indicated by saying that  the Schauder frame 
$\{(x_j, x_j^*)\}_{j=1}^\infty$ is \emph{$K$-unconditional}.

	 \subsection{}
We will use $\{r_n(t)\}_{n=1}^{\infty}$ to denote the sequence of Rademacher functions 
defined on the segment $[0,1]$ by $r_n(t) = \sign \sin (2^n \pi t)$.
The next lemma formulates the fact that for $1\le p\le 2$ the space $L^p$ has cotype $2$.
	
	\begin{lemma}
			\label{cotype}
	Let $1\le p \le 2$.	There exists a constant $C_p$ such that for any finite 
	number of functions $f_1, f_2, \ldots f_N \in L^p(\R)$ we have 
		\begin{equation}
		\label{eq:cotype}
		\Big( \sum_{j=1}^N \|f_j\|_p^2 \Big)^{1/2}\leq C_p \int_0^1 \Big\|\sum_{j=1}^N r_j(t) f_j\Big\|_p\, dt.
		\end{equation}
	\end{lemma}
	
This is a key point where the assumption   $1 \le p \le 2$ is used
in the proof. Note that inequality \eqref{eq:cotype} does not hold for $p>2$.
	
	A proof of \lemref{cotype}, as well as other information about the 
	notions of type and cotype of a Banach space,
	can be found e.g.\ in \cite[Chapter III.A]{Wojt}. 
	In fact, $C_p$ may be chosen to be an absolute constant not depending on $p$, but
	 we will not use this fact.

	 \subsection{}	
	A set $\Lambda \subset \R$ is said to be \emph{uniformly discrete} if there exists $\delta>0$ such that 
	$|\lam' - \lam| \geq \del$ for any two distinct points $\lam, \lam' \in \Lam$. 
	We will use the following elementary fact (for a proof, see \cite[Proposition 2.1]{OSSZ}).
	
	\begin{lemma}
		 	\label{sum_of_norms}
		Let  $f \in L^p(\R)$, $1 \le p < \infty$, and suppose that
		 $\Lambda\subset \R$ is  a 
		  uniformly discrete set. Then for any bounded interval $I\subset\R$ we have
		 	\begin{equation}
		\sum_{\lambda\in\Lambda} \int_I |T_\lambda f|^p <+\infty.
		\end{equation}
	\end{lemma}


	\section{Proof of \thmref{thmA1}}	
		\label{secP2}

	Let $1 \le p \le 2$, and assume that there exists
in the space $L^p(\R)$ an unconditional Schauder frame of the form
$\{(T_{\lam_j} g, g_j^*)\}_{j=1}^{\infty}$  where $g \in L^p(\R)$,
$\{\lam_j\}$ are real numbers, and
$\{g_j^*\}$ are continuous linear functionals on $L^p(\R)$.
Our goal is to show that this assumption leads to a contradiction.
 With no loss of generality we may assume that $\|g\|_p=1$.

	A standard approach involves an application of \lemref{unc} which implies that
	\begin{equation}
	\label{eqA3}
	\Big\|  \sum_{j=1}^{N} r_j(t) g_j^*(f) T_{\lambda_j} g \Big\|_p \leq K \|f\|_p,
	\quad t \in [0,1],
	\end{equation}
	where $\{r_j(t)\}$ is the sequence of  Rademacher functions. If we now take 
	$\int_0^1 dt$ and use \lemref{cotype}, 	and finally also take the limit as $N \to \infty$, 
	we obtain	the $\ell^2$-estimate
	\begin{equation}
	\label{eststd}
		\Big( \sum_{j=1}^{\infty} |g_j^*(f)|^2  \Big)^{1/2} \leq C_p K \|f\|_p
	\end{equation}	
	(note that we have used the assumption  $\|g\|_p=1$).
	However this estimate for  the coefficients $\{g_j^*(f)\}$ is not good 
	enough to yield the full statement of \thmref{thmA1},
	and so our approach will involve a certain improvement of the inequality \eqref{eststd}.

	We will present the proof in several steps.
	
	\subsection{Step 1: An $\ell^1$-estimate for blocks of coefficients}
		Let us fix $\delta > 0$ small enough 
		such that   $\|T_a g - T_b g\|_p \leq 1/2$
	whenever $a,b$ are two real numbers satisfying $|a-b|\le \delta$. 
	We consider a partition of $\R$ into
	intervals $ [k\delta, (k+1)\delta)$, $k \in \Z$.  It induces a partition of the
	unconditional Schauder frame into ``blocks'' 
	$\{(T_{\lambda_j} g, g_j^*)\}_{j \in A_k}$, where 
	\begin{equation}
	A_k := \{j: \lambda_j\in [k\delta, (k+1)\delta)\}, \quad k \in \Z.
	\end{equation}
It is important to keep in mind that some of the sets $A_k$ may be infinite, since 
the sequence $\{\lam_j\}_{j=1}^{\infty}$ does not necessarily tend to $\pm \infty$
and thus a bounded interval may in principle contain infinitely many elements $\lam_j$ of the sequence.

	For any $f\in L^p(\R)$ we thus have the unconditionally convergent series expansion 
	\begin{equation}
	\label{represent}
	f=\sum_{k\in\Z}\sum_{j\in A_k} g_j^*(f) T_{\lambda_j}g.
	\end{equation}

	We now make the following observation: given arbitrary nonnegative scalars $\{a_j\}_{j\in A_k}$ such that
	 only finitely many of them are nonzero, we have
	\begin{equation}
	\label{key0}
	\Big\|\sum_{j\in A_k} a_j T_{\lambda_j}g\Big\|_p\ge \frac{1}{2}\sum_{j\in A_k}a_j.
	\end{equation}	
	Indeed, using the fact that $\|T_{k\delta}g - T_{\lambda_j}g \|_p \le 1/2$
	for each $j\in A_k$, we obtain
	\begin{align}
	\Big\|\sum_{j\in A_k} a_j T_{\lambda_j}g\Big\|_p &\ge \Big\| \sum_{j\in A_k} a_j T_{k\delta} g \Big\|_p - \Big\| \sum_{j\in A_k} a_j (T_{k\delta} g -T_{\lambda_j}g)\Big\|_p\\ 
	&=\Big(\sum_{j\in A_k} a_j\Big) - \Big\| \sum_{j\in A_k} a_j (T_{k\delta} g -T_{\lambda_j}g)\Big\|_p \ge\frac{1}{2}\sum_{j\in A_k} a_j
\end{align}
	(again using $\|g\|_p=1$) which shows that \eqref{key0} is valid.
	
	Assume that for some $k \in \Z$ the set $A_k$ is finite. Then, if
	 we put $a_j=|g_j^*(f)|$ in the inequality \eqref{key0}, we arrive at the estimate
	\begin{equation}
		\label{simple_useful}
			\sum_{j\in A_k} |g_j^*(f)| \leq 2 \, \Big\| \sum_{j\in A_k} |g_j^*(f)| \; T_{\lambda_j} g \Big\|_p.
	\end{equation}
	
	Next we observe that \eqref{simple_useful}
	holds also if $A_k$ is an infinite set. Indeed, the series
	on the right hand side is guaranteed to 
	converge unconditionally due to \lemref{unc}, 
	so the inequality can be established 
	by applying  the estimate \eqref{simple_useful} to 
	the coefficients $\{g_j^*(f)\}$ in increasing finite subsets
	of $A_k$, and passing to the limit.

	Moreover, again by \lemref{unc} we have
	\begin{equation}
	\label{simple_est}
		\Big\| \sum_{j\in A_k} |g_j^*(f)| \; T_{\lambda_j} g \Big\|_p \leq K \|f\|_p,
	\end{equation}
	hence \eqref{simple_useful} and \eqref{simple_est} imply the $\ell^1$-estimate 
	\begin{equation}
		\label{one_block}
		\sum_{j\in A_k} |g_j^*(f)| \leq 2 K \|f\|_p
	\end{equation}
	for the coefficients in the $A_k$ block of the series expansion  \eqref{represent}.

	\subsection{Step 2: Compactness of   partial sum operators}
For each $N$ we let
\begin{equation}
\label{eqC1}
S_N f = \sum_{|k|\le N} \sum_{j\in A_k} g_j^*(f) T_{\lambda_j} g, \quad f \in L^p(\R),
\end{equation}
be the $N$'th partial sum of the series expansion \eqref{represent}. 
We note that some of the blocks  $A_k$ may be infinite, hence the
series \eqref{eqC1} may contain infinitely many terms. However, the series
converges unconditionally due to \lemref{unc}.

We consider $S_N$ as a linear operator on the space $L^p(\R)$. 
Our goal  is to show that $S_N$ 
is a compact operator. This is obvious if all the blocks $A_k$ happen to 
be finite, since in this case   $S_N$ is an operator  of finite rank.
In the general case, where some of the blocks $A_k$ may be infinite, 
we argue as follows.

It would be enough to prove that for each fixed $k \in \Z$, the operator  $B_k$ defined by
\begin{equation}
\label{eqC2}
B_k f = \sum_{j\in A_k} g_j^*(f) T_{\lambda_j} g
\end{equation}
is compact. We will 
show that $B_k$ can be approximated by  operators of finite rank.
Let $\eps > 0$, and choose $\eta > 0$ small enough 
so that we have  $\|T_a g - T_b g\|_p \le \eps$
	whenever $a,b$ are two real numbers satisfying $|a-b|\le \eta$.
	Let us take a finite $\eta$-net for the segment $[k\delta, (k+1)\delta)$, 
	that is, a finite set $F$ such that each point of 
	the segment $[k\delta, (k+1)\delta)$ lies within distance at most $\eta$
	from some point of $F$. For each $j\in A_k$ we choose 
	a point $\sig_j$ in this $\eta$-net $F$ satisfying 
	$|\sig_j - \lambda_j |\le \eta$. Then
$  \|T_{\sig_j} g - T_{\lambda_j} g \|_p \leq \eps$, for each $j \in A_k$.

		Define an operator $E_k$ by
	\begin{equation}
	\label{def_rop}
	E_k f = \sum_{j\in A_k} g_j^*(f) (T_{\sig_j} g - T_{\lambda_j} g), \quad f \in L^p(\R).
	\end{equation}
	We show that the series \eqref{def_rop} converges unconditionally and defines a continuous linear operator 
	$E_k : L^p(\R) \to L^p(\R)$. Indeed,  using the estimate \eqref{one_block} we obtain
	\begin{equation}
		\|E_k f\|_p \le 
			 \sum_{j\in A_k} |g_j^*(f)| \cdot \|T_{\sig_j} g - T_{\lambda_j} g \|_p 
		\leq \eps  \sum_{j\in A_k} |g_j^*(f)| \le 2K \eps \|f\|_p,
	\end{equation}
which  establishes the unconditional convergence of 
the series \eqref{def_rop} as well as the continuity of the operator $E_k$.
Moreover, we obtain that $\|E_k\| \le 2K\eps$.

Note that we have
	\begin{equation}
	\label{exp_sop}
	(B_k + E_k) f =   \sum_{j\in A_k}   g_j^*(f) T_{\sig_j} g, \quad f \in L^p(\R),
	\end{equation}
	the series \eqref{exp_sop} again being unconditionally convergent. 
But recall that each  $\sig_j$ belongs to 
the finite $\eta$-net $F$ of the segment $[k\delta, (k+1)\delta)$.
Hence the image of the operator $B_k + E_k$ lies in the 
finite-dimensional linear subspace spanned by the translates
$\{T_{\lam} g\}_{\lam \in F}$, which implies that
$B_k + E_k$ is an operator of finite rank.

	We have thus shown that for any $\eps > 0$ there exists an operator  $B_k + E_k$
of finite rank,  that approximates $B_k$
in the operator norm with error at most $2K\eps$.
 Hence $B_k$ is a compact  operator, for each $k \in \Z$.
As a consequence, $S_N$ must be compact.

	\subsection{Step 3: An application of cotype $2$ inequality}
	
	Our next goal is to use the estimate \eqref{simple_useful} 
	 in order to obtain an improvement of  the inequality   \eqref{eststd}. 
	 The key observation is that we will get a better estimate if instead of multiplying each 
	 summand by a random sign separately as in \eqref{eqA3},
	  we will do it for each block. 
	  
	Indeed, by an application of \lemref{unc} we have
	\begin{equation}
	\label{eqB8}
		\Big\| \sum_{|k| \leq N} r_{n_k}(t) \sum_{j\in A_k} |g_j^*(f)| \; T_{\lambda_j} g \Big\|_p \leq K \|f\|_p,
	\quad t \in [0,1],
	\end{equation}
	where $\{r_{n_k}(t)\}_{k \in \Z}$ is some  enumeration of the
	sequence of  Rademacher functions. If we take 
	$\int_0^1 dt$ and use \lemref{cotype}, we obtain
	\begin{equation}
	C_p^{-1}\Big( \sum_{|k| \leq N} \Big\| \sum_{j\in A_k} |g_j^*(f)| \; T_{\lambda_j} g 
	\Big\|_p^2 \Big)^{1/2}\leq K \|f\|_p.
	\end{equation}
	Finally, combining this inequality together with \eqref{simple_useful} and 
	taking the limit as $N \to \infty$, we arrive at the key estimate
	\begin{equation}
		\Big[ \sum_{k \in \Z} \Big( \sum_{j\in A_k} |g_j^*(f)| \Big)^2 \Big]^{1/2}\leq 2K C_p \|f\|_p,
		\label{key}
	\end{equation}
	which is a ``mixed'' $(\ell^1,\ell^2)$-estimate for  the coefficients $\{g_j^*(f)\}$ that
	improves  on  the standard $\ell^2$-estimate given in \eqref{eststd}.

\subsection{Step 4: Compactness of restriction operator}

Let $I\subset\R$ be an arbitrary bounded interval, and let 
$R_I$ be the operator of 
restriction to the interval $I$, that is, $R_I f := f|_I$. Then 
$R_I: L^p(\R) \to L^p(I)$ is a continuous linear operator.

Consider again the sequence of partial sum operators $\{S_N\}$ defined by
\eqref{eqC1}. Then  for every  $f \in L^p(\R)$ we
have  $(S_N f)|_I \to f|_I$ in $L^p(I)$, that is,
$R_I \circ S_N \to R_I $  in the strong operator topology.
If we show that $\{R_I \circ S_N\}$ is a Cauchy sequence of operators  $L^p(\R) \to L^p(I)$
 (with respect to the operator norm), then it would follow that $R_I \circ S_N \to R_I$ 
in the operator norm. But  as we have shown above, 
 $S_N$,  and hence also $R_I \circ S_N$, are compact operators.
We will therefore obtain that  $R_I$ must be  a compact operator as well, 
which leads to a contradiction.

(We note that this approach,
based on compactness of the restriction operator $R_I$,
goes back to \cite{OSSZ} and has been used also in  \cite{FOSZ}, \cite{BerCar}.)

It remains to prove that 
$\{R_I \circ S_N\}$ is a Cauchy sequence. For any $M>N$ we have
\begin{align}
	& \|R_I \circ (S_M - S_N) f\|_p = \Big\| \sum_{N<|k|\le M} \sum_{j\in A_k} g_j^*(f) (T_{\lambda_j}g)|_I \Big\|_p \\
	& \qquad  \qquad \leq \sum_{N<|k|\le M} \sum_{j\in A_k} |g_j^*(f)| \cdot \|(T_{\lambda_j}g)|_I\|_p \\
	& \qquad \qquad  \leq \sum_{N<|k|\le M}  \Big( \sum_{j\in A_k} |g_j^*(f)| \Big) \cdot \max_{ \lam \in [k\delta, (k+1)\delta]} \|(T_{\lam} g)|_I\|_p. \label{est1}
\end{align}

For each  $k\in\Z$ we can choose a point $\lambda^{(k)}$ in the segment
$[k\delta, (k+1)\delta]$  such that 
\begin{equation}
\label{eqB5}
\max_{\lam \in [k\delta, (k+1)\delta]}\|(T_{\lambda}g)|_I\|_p = \|(T_{\lambda^{(k)}} g)|_I\|_p.
\end{equation}

Assume first that $1 < p \leq 2$, and let $p' = p/(p-1)$ be the exponent  conjugate to $p$.
Using \eqref{eqB5} and H\"{o}lder's inequality we obtain that the quantity \eqref{est1} is bounded by 
\begin{align}
	& \sum_{N<|k|\le M} \Big( \sum_{j\in A_k} |g_j^*(f)| \Big) \|(T_{\lambda^{(k)}} g)|_I\|_p\\
	&\qquad \leq \Big[ \sum_{k\in\Z} \Big( \sum_{j\in A_k} |g_j^*(f)| \Big)^{p'} \Big]^{1/p'} \Big[ \sum_{N<|k|\le M} \int_I |T_{\lambda^{(k)}}g|^p \Big]^{1/p}. \label{est2}
\end{align}
We observe that the sequence  $\{\lambda^{(k)}\}_{k\in\Z}$ can be partitioned into 
two uniformly discrete subsequences (since each one of $\{\lambda^{(2k)}\}$ and $\{\lambda^{(2k+1)}\}$ is uniformly discrete). 
Hence \lemref{sum_of_norms} implies that for any $\eps>0$ there exists $N(\eps)$ 
such that the second factor  in \eqref{est2}
will be smaller than $\eps$ for $M > N > N(\eps)$. 
The first factor can be estimated using the fact that $p'\ge 2$ and inequality \eqref{key}. 
It implies that  for $M>N > N(\eps)$ the quantity \eqref{est2} does not exceed
\begin{equation}
 \Big[ \sum_{k\in\Z} \Big( \sum_{j\in A_k} |g_j^*(f)| \Big)^2 \Big]^{1/2}  \eps \le 2KC_p\,\eps\|f\|_p.
\end{equation}
Therefore, for $M>N>N(\eps)$ we have 
$\|R_I \circ (S_M-S_N)\| \le 2KC_p \, \eps$,
which shows that $\{R_I \circ S_N\}$ is indeed a Cauchy sequence
and leads to the desired contradiction.

	Lastly, in the remaining case where $p=1$, we   simply replace
	 the first factor in \eqref{est2} with
	$\sup_{k\in\Z} \big\{ \sum_{j\in A_k} |g_j^*(f)|  \big\}$
	and the proof continues in the same way. However we note
	that in this case a more general result is actually true: it is known that
	in the space $L^1(\R)$ there do not exist any unconditional Schauder frames, 
	whether formed by translates or not (see \cite[Section 4.3]{BerCar}).


	\section{Remarks}
		\label{secR1}

	\subsection{}
	\thmref{thmA1} can be easily generalized to the space $L^p(\R^d)$. 
	It is also easy to extend the result to unconditional Schauder frames formed 
	by translates of a finite number of functions. That is, for $1 \le p\le 2$ there does not exist 
in the space $L^p(\R^d)$ any unconditional Schauder frame of the form
$\{(g_j, g_j^*)\}$ such that each $g_j$ is a translate of a function
from some finite set $G \sbt L^p(\R^d)$.

	\subsection{}	
	In fact our proof establishes the following more general
	version of \thmref{thmA1}.
	
		\begin{thm}
	\label{thmA8}
	Assume that the system $\{(T_{\lambda_j}g, g_j^*)\}_{j=1}^{\infty}$ 
	forms an unconditional Schauder frame for a (closed, linear) 
	subspace $X \sbt L^p(\R)$, $1 \le p \le 2$.	
	 Then for every bounded interval $I \sbt \R$ the restriction operator 
	 $R_I: X \to L^p(I)$ is compact.
	\end{thm}
	
	Moreover, in this case the subspace $X$ embeds into $\ell^p$, see \cite[Proposition 5.3]{FOSZ}.

	\subsection{} \label{seqC3}
		Let $q \ge 2$, and let $p = q/(q-1)$ be the exponent  conjugate to $q$.	
	By identifying the dual space of $L^q(\R)$ with $L^p(\R)$ in the usual way,
	one can define the translates $T_{\lam} g^*$, $\lam \in \R$, of 
	a continuous linear functional $g^*$ on $L^q(\R)$. 
	This allows to consider Schauder frames in $L^q(\R)$ of the form
	$\{(g_n, T_{\lam_n} g^*)\}$, namely, where the coefficient functionals
	are translates of a single functional $g^*$.
	As a consequence of \thmref{thmA1} we can then obtain the following result
	(compare with \cite[Corollary 3.3]{FOSZ}):  
	
		\begin{thm}
	\label{thmA2}
		For   $q \ge 2$ there does not exist 
in the space $L^q(\R)$ any unconditional Schauder frame of the form
$\{(g_n, T_{\lam_n} g^*)\}$  where $g_n \in L^q(\R)$,
$g^*$ is a continuous linear functional on $L^q(\R)$, and 
$\{\lam_n\}$ are real numbers.
	\end{thm}

This can be deduced from \thmref{thmA1} using the following 
fact (see, for example, the discussion before the proof of \cite[Corollary 3.3]{FOSZ}):  if 
$\{(x_j, x_j^*)\}_{j=1}^{\infty}$ is an unconditional Schauder frame in 
		 a \emph{reflexive} Banach space $X$, then the system 
		 $\{(x_j^*, x_j)\}_{j=1}^{\infty}$ forms an
		 unconditional Schauder frame for the dual space $X^*$.
		 
We outline a short proof of the latter fact. It is straightforward to verify that every
$x^* \in X^*$ has a weakly unconditionally convergent expansion
$x^* = \sum_{j=1}^\infty x^*(x_j) x_j^*$
(due to reflexivity, the weak convergence in $X^*$ is the same as the weak-$*$ convergence). On the other hand, $X^*$ being reflexive does not contain any subspace isomorphic to $c_0$
(see \cite[Theorem II.A.14]{Wojt}), hence every weakly unconditionally 
convergent series in $X^*$ is unconditionally convergent 
(see \cite[Proposition II.D.5]{Wojt}).

\subsection{}
One can consider also translates of a function $g$ on the circle $\T = \R / \Z$.
By an argument similar to \cite[Theorem 1]{OlsZal} one can show that
the space $L^p(\T)$, $1 \le p < \infty$, does not admit
a Schauder basis formed by translates of a single function.
On the other hand, it can be deduced from
\cite[Section 4]{FPT} that there do exist Schauder frames
of translates in $L^p(\T)$ for any $1\le p < \infty$.
Our technique allows us to establish that Schauder
frames of this form cannot be unconditional:

\begin{thm}
There does not exist in the space
  $L^p(\T)$, $1 \le p < \infty$,
 any unconditional Schauder frame 
of translates, i.e.\ of the form
$\{(T_{\lam_n} g, g_n^*)\}$  where $g \in L^p(\T)$,
$\{\lam_n\} \subset \T$,
and $\{g_n^*\}$ are continuous linear functionals on $L^p(\T)$.
\end{thm}

This can be proved by adapting Steps 1 and 2 from
the proof of \thmref{thmA1}. Based on these steps
one can show that the existence of   an unconditional Schauder frame 
of translates
 in $L^p(\T)$ implies compactness of the identity operator, a contradiction.

\subsection{}
A system of vectors $\{x_n\}$ in a Hilbert space $H$ is called
a \emph{frame} if there are positive constants $A,B$ such that 
the inequalities 
$A \|x\|^2 \le \sum_n |\dotprod{x}{x_n}|^2 \leq B\|x\|^2$
hold for every $x \in H$.
In this case there exists another frame $\{y_n\}$ (called the ``dual''  frame)
such that the system 
$\{(x_n, y_n)\}$ forms an unconditional Schauder frame in $H$,
where the dual space $H^*$ is  identified with $H$ in the usual way
(see \cite[Section 4.7]{You01}).

It was proved in \cite{CDH} that
in the space $L^2(\R)$ there does not exist any frame formed
by translates of a single function. This result can now be viewed
as a special case of \thmref{thmA1}.

The following open problem was communicated to us by D. Freeman:
If $\{(x_n, y_n)\}$ is an unconditional Schauder frame in a Hilbert space $H$,
does there exist a sequence of nonzero scalars $\{\alpha_n\}$ such that
 $\{\alpha_n x_n\}$ and $\{ \alpha_n^{-1} y_n\}$ are both frames in $H$?
For more information related to this problem we refer the reader to  \cite{BFPS24}.

In this direction, the following result was obtained in \cite[Theorem 3.13]{HLLL14}.
 
 		\begin{prop}
	\label{propC1}
	Let $\{(x_n, y_n)\}$ be an unconditional Schauder frame in a  Hilbert space $H$. 
	If both $\inf_n  \|x_n\|$ and $\inf_n  \|y_n\|$ are nonzero, 
		  then $\{x_n\}$ and $\{y_n\}$ are both frames.
	\end{prop}

We use  this opportunity to present a short, self-contained proof of this fact.

\begin{proof}
Due to \lemref{unc} we have
 $	\|\sum_{n} r_n(t) \dotprod{x}{y_n} x_n \|^2 \leq K^2 \|x\|^2$,
 $x \in H$, 
 	where $\{r_n(t)\}$ is the sequence of  Rademacher functions. If we take 
	$\int_0^1 dt$ then orthogonality of Rademacher functions gives 
	$\sum_{n} |\dotprod{x}{y_n}|^2 \|x_n\|^2 \leq K^2 \|x\|^2$.
Using the assumption that $\inf_n  \|x_n\| > 0$ we conclude
that $\{y_n\}$ is a Bessel sequence (see \cite[Section 4.2]{You01}).

The system $\{(y_n, x_n)\}$ is also 
an unconditional Schauder frame in $H$ (see \secref{seqC3}),
hence by symmetry also $\{x_n\}$ is a Bessel sequence.
We can therefore define a  continuous  linear operator $A: H\to \ell^2$ 
by $Ax=\{\dotprod{x}{x_n}\}$. The
 adjoint operator $A^*: \ell^2 \to H$ is given by
$A^*(\{c_n\}) = \sum_{n} c_n x_n$.
The fact that $\{(x_n, y_n)\}$ is a Schauder frame 
and $\{y_n\}$ is a Bessel sequence implies
that  $A^*$ is surjective, 
which in turn is equivalent to the existence of a constant $c>0$
 such that $\|A x\|^2 \ge c\|x\|^2$, $x \in H$  (see \cite[I.A.13]{Wojt}).
 	This inequality together with the boundedness of $A$
 	means that $\{x_n\}$ is a frame. Finally, by symmetry also
 	$\{y_n\}$ is a frame, and so the claim is proved.
\end{proof}


\end{document}